\documentclass[11pt,oneside]{amsart}

\usepackage{amsfonts}
\usepackage{amsmath, amsthm, amssymb, ulem, amscd}

\def\crn#1#2{{\vcenter{\vbox{
        \hbox{\kern#2pt \vrule width.#2pt height#1pt
           }
          \hrule height.#2pt}}}}
\def\intprod{\mathchoice\crn54\crn54\crn{3.75}3\crn{2.5}2}
\def\into{\mathbin{\intprod}}

\newcommand{\stopthm}{\hfill$\square$\medskip}

\newcommand{\Ric}{\operatorname{Ric}}
\newcommand{\End}{\operatorname{End}}
\newcommand{\Aut}{\operatorname{Aut}}

\newcommand{\Hol}{\operatorname{Hol}}

\newcommand{\ev}{\operatorname{ev}}
\renewcommand{\span}{\operatorname{span}}

\newcommand{\contr}{\operatorname{contr}}

\newcommand{\R}{\mathbb R}
\newcommand{\N}{\mathbb N}
\newcommand{\C}{\mathbb C}

\newcommand{\Rt}{\widetilde{R}}

\newcommand{\Vt}{\widetilde{V}}
\newcommand{\nt}{\widetilde{\nabla}}

\newcommand{\xt}{\widetilde{\xi}}
\newcommand{\zt}{\widetilde{\zeta}}
\newcommand{\et}{\widetilde{\eta}}

\newcommand{\Dt}{\widetilde{\Delta}}

\newcommand{\eb}{\overline{\eta}}

\newcommand{\gt}{\widetilde{g}}

\newcommand{\cE}{\mathcal{E}}

\newcommand{\cG}{\mathcal{G}}
\newcommand{\cGt}{\widetilde{\mathcal{G}}} 
\newcommand{\cEt}{\widetilde{\mathcal{E}}} 
\newcommand{\cRt}{\widetilde{\mathcal{R}}}

\newcommand{\cV}{\mathcal{V}}
\newcommand{\cT}{\mathcal{T}}
\newcommand{\cD}{\mathcal{D}}

\newcommand{\fX}{\mathfrak{X}}
\newcommand{\fh}{\mathfrak{h}\mathfrak{o}\mathfrak{l}}
\newcommand{\fht}{\widetilde{\mathfrak{h}\mathfrak{o}\mathfrak{l}}}

\theoremstyle{plain}
\newtheorem{theorem}{Theorem}[section]

\newtheorem{corollary}[theorem]{Corollary}

\theoremstyle{definition}

\theoremstyle{remark}

\numberwithin{equation}{section}

\title[Conformal Holonomy Equals Ambient Holonomy]{Conformal Holonomy
  Equals Ambient Holonomy}     

\author{Andreas \v{C}ap} 
\address{Faculty of Mathematics\\
University of Vienna\\
Oskar--Morgenstern--Platz 1\\
1090 Wien, Austria}
\email{Andreas.Cap@univie.ac.at}  

\author{A. Rod Gover}
\address{The University of Auckland\\
Private Bag 92019\\
Auckland 1142, New Zealand}

\email{r.gover@auckland.ac.nz}

\author{C. Robin Graham}
\address{Department of Mathematics, University of Washington,
Box 354350\\
Seattle, WA 98195-4350, USA}
\email{robin@math.washington.edu}

\author{Matthias Hammerl}
\address{Department of Mathematics and Informatics\\
University of Greifswald\\
Walther-Rathenau-Str. 47\\
17489 Greifswald, Germany}
\email{matthias.hammerl@uni-greifswald.de}

\begin{document}

\begin{abstract}
This paper studies the relation between two notions of holonomy on a
conformal manifold.  The first is the conformal holonomy, defined to be the
holonomy of the normal tractor connection.  The second is the holonomy of
the Fefferman-Graham ambient metric of the conformal manifold.  It is shown
that the infinitesimal conformal holonomy and the infinitesimal ambient
holonomy always agree up to the order that the ambient metric is defined.  
\end{abstract}

\maketitle

\thispagestyle{empty}

\renewcommand{\thefootnote}{}
\footnotetext{ A\v C \& ARG  gratefully acknowledge support from the Royal
  Society of New Zealand via Marsden Grant 13-UOA-018.  
Further, support by projects P23244-N13 (A\v C and MH) and P27072-N25 
(A\v C) of the Austrian Science Fund (FWF) is gratefully acknowledged.   
Research of CRG is partially supported by NSF grant \# DMS 1308266.}    
\renewcommand{\thefootnote}{1}

\section{Introduction}\label{intro}

The tractor bundle $\cT$ of a smooth conformal manifold $(M,c)$ of
dimension $n\geq 3$ and signature $(p,q)$, $p+q=n$, is a rank $n+2$
vector bundle naturally associated to the conformal structure, which
carries a canonical connection $\nabla$, see \cite{BEG}.  This
connection is characterized by a normalization condition on its
curvature, whence it is called the normal tractor connection 
\cite{CG}. It can be viewed as a conformally invariant analog of the
Levi-Civita connection of a Riemannian manifold and has played an
essential role in many recent developments in conformal geometry.  The
holonomy of $(\cT,\nabla)$ is called the conformal holonomy of
$(M,c)$; following early work including classification results
\cite{Armstrong}, \cite{Leitner}, its study has been the focus of active
recent research.  

Another invariant object associated to a conformal manifold is the
ambient metric of \cite{FG1}, \cite{FG2}.  This is a smooth
pseudo-Riemannian metric of signature $(p+1,q+1)$ on a space of dimension
$n+2$, determined up to 
diffeomorphism along a canonical hypersurface, to infinite order if $n$ is  
odd, and to order $n/2-1$ if $n$ is even.  Its Levi-Civita connection is 
another connection associated to the conformal manifold and one can also
consider its holonomy.  Because the holonomy group of a connection is a
global invariant and the  
ambient metric is only invariantly defined as a jet along a hypersurface,
its holonomy group is not the appropriate object to study.  Instead we  
consider the infinitesimal holonomy, which depends only on the jet at a   
point. 
The main result of this paper asserts that, suitably interpreted, the
infinitesimal holonomies of the tractor connection and the Levi-Civita 
connection of the ambient metric agree at each point.  

In order to formulate the result precisely, we describe a realization 
of the tractor bundle in ambient terms which was derived in \cite{CG}.
Details will be provided in \S\ref{ambtrac}.  
If $(M,c)$ is a conformal manifold, its metric bundle is the ray bundle 
$\cG\subset S^2T^*M$ whose sections are the metrics $g\in c$.  The ambient
space is $\cG\times \R$, in which $\cG$ is embedded as the hypersurface
$\cG\times \{0\}$.  There are dilations $\delta_s:\cG\rightarrow \cG$ given
by $\delta_s(x,g_x)=(x,s^2g_x)$, $s>0$, which extend to $\cG\times \R$
acting in the first factor.  For $x\in M$, we denote by $\cG_x$ the fiber
of $\cG$ over $x$, and we view $\cG_x$ as a 1-dimensional submanifold of
$\cG\times \R$ via $\cG_x\subset \cG = \cG\times \{0\}\subset
\cG\times \R$.  Then $T(\cG\times \R)\,\big{|}_{\cG_x}$ denotes the tangent
bundle to $\cG\times \R$ restricted to the submanifold $\cG_x$, a rank
$n+2$ vector bundle over $\cG_x$.  The standard tractor bundle of $(M,c)$
can be realized as the rank $n+2$ vector bundle  
$\cT\rightarrow M$ with fiber      
\begin{equation}\label{tracdef}
\cT_x=\left\{U\in\Gamma\big(T(\cG\times \R)\,\big{|}_{\cG_x}\big): 
(\delta_s)^*U=s^{-1}U,\;s>0\right\}.     
\end{equation}
Thus a section of $\cT$ on $M$ is a vector field in $\cG\times\R$ defined
on $\cG$ which is homogeneous of degree $-1$ with respect to the
$\delta_s$.  

As we will also review in \S\ref{ambtrac}, an ambient metric for $(M,c)$ 
is a pseudo-Riemannian metric $\gt$ which 
is defined in a dilation-invariant neighborhood $\cGt$ of $\cG$ in
$\cG\times \R$ by certain conditions.  As indicated above, it is uniquely
determined by $(M,c)$ up to diffeomorphism to infinite order if $n$ is  
odd and to order $n/2-1$ if $n$ is even.  

We recall the notion of infinitesimal holonomy.  
A standard reference is \cite{KN}, where however the formulation is in
terms of a connection on a principal bundle instead of a connection on an
associated vector bundle as considered here.   
If $(\cV,\nabla)$ is a smooth vector bundle with  
connection on a manifold $M$ and $x\in M$, the   
infinitesimal holonomy algebra $\fh_x$ of $(\cV,\nabla)$ at $x$ is the
subspace of $\End\cV_x$ defined by 
\begin{equation}\label{holdef}
\fh_x=\span_\R\left\{\nabla_{\eta_k}\nabla_{\eta_{k-1}}\cdots
\nabla_{\eta_3}\big(R(\eta_1,\eta_2)\big)(x): k\geq 2
,\,\,\eta_1,\ldots, \eta_k\in \fX(M)\right\}.
\end{equation}
Here $\fX(M)$ denotes the space of smooth vector fields on $M$ 
and $R:\Lambda^2TM\rightarrow \End\cV$ the curvature 
of $\nabla$.  It is a standard fact that $\fh_x$ is a subalgebra of
$\End\cV_x$ for its natural Lie algebra structure with bracket the
commutator of endomorphisms.  Clearly $\fh_x$ depends only on the infinite
order jet of $\nabla$ at $x$, and so in particular there is generally no 
relation between $\fh_x$ and $\fh_y$ for $x\neq y$.  However, if $M$ and
$(\cV,\nabla)$ are real-analytic, then  
$\fh_x$ is the Lie algebra of $\Hol_x$, where $\Hol_x\subset \Aut \cV_x$
is the usual 
holonomy group of $(\cV,\nabla)$ defined by parallel translation around
loops based at $x$.  Of course, $\Hol_x$ is always isomorphic to $\Hol_y$
for $M$ smooth and connected.  

For a conformal manifold $(M,c)$, we denote by $\fh_x$ the 
infinitesimal holonomy at $x$ of $(\cT,\nabla)$, where $\nabla$ is the
normal tractor connection.  Thus $\fh_x$ is a subalgebra of $\End \cT_x$.
The realization \eqref{tracdef} of $\cT_x$ induces the realization  
\begin{equation}\label{endos}
\End \cT_x=\left\{E\in\Gamma\big(\End T(\cG\times 
\R)\,\big{|}_{\cG_x}\big): (\delta_s)^*E=E,\;s>0\right\}
\end{equation}
of $\End \cT_x$.  Thus 
an element of $\fh_x$ is realized as a section of the vector bundle $\End  
T(\cG\times \R)|_{\cG_x}$ over $\cG_x$ which is homogeneous of degree 0 with  
respect to the $\delta_s$.  For any $z\in \cG_x$, evaluation at $z$ is an
isomorphism 
$$
\ev_z:\End\cT_x\rightarrow \End T_z(\cG\times\R).
$$  
So $\ev_z(\fh_x)$ is an isomorphic copy of $\fh_x$ in 
$\End T_z(\cG\times\R)$.         

If $\gt$ is an ambient metric for $(M,c)$ and $x\in M$, the infinitesimal
holonomy at $z\in \cG_x$ of the Levi-Civita connection $\nt$ of  
$\gt$ is a subalgebra of $\End T_z\cGt=\End T_z(\cG\times \R)$.  
If $n$ is odd, we denote this subalgebra $\fht_z$.  This is clearly
independent of the infinite-order ambiguity in $\gt$.  However, when $n$ is   
even, the ambient metric is determined by $(M,c)$ only to order $n/2-1$  
along $\cG$.  So we need to restrict the number of 
differentiations transverse to $\cG$ to avoid 
this ambiguity.  Therefore, when $n\geq 4$ is even, we define 
\begin{equation}\label{ambhol}
\fht_z=\span_\R\left\{\nt_{\xt_k}\nt_{\xt_{k-1}}\cdots
\nt_{\xt_3}\left(\Rt(\xt_1,\xt_2)\right)(z):k\geq
2,\,\,\xt_1,\ldots,\xt_k\in \fX(\cGt)\right\},  
\end{equation}
where $\Rt$ is the curvature of $\nt$, but we
impose the requirement that no more than $n/2-2$ of  
the vector fields $\xt_1,\ldots,\xt_k$ are somewhere transverse to $\cG$.
Then $\nt_{\xt_k}\nt_{\xt_{k-1}}\cdots  
\nt_{\xt_3}\left(\Rt(\xt_1,\xt_2)\right)$ depends on at most $n/2-1$    
transverse derivatives of $\gt$, so its value at $z$ is independent of the
ambiguity at order $n/2$.   

Our main result is the following.  
\begin{theorem}\label{main}
Let $(M,c)$ be a conformal manifold of dimension $n\geq 3$ and $\gt$ an 
ambient metric for $(M,c)$.  If $x\in M$ and $z\in \cG_x$, then
$$
\ev_z(\fh_x) = \fht_z.   
$$
\end{theorem}

An immediate corollary is the equality of restricted tractor and ambient
holonomy groups in the 
odd-dimensional real-analytic case.  Recall that if $(\cV,\nabla)$ is a
vector bundle with connection on a smooth manifold $M$ and $x\in M$, then
the restricted holonomy group is 
$$
\Hol_x^0(\cV,\nabla)=\{L_\gamma\}\subset \Aut \cV_x, 
$$
where $\gamma$ is a smooth contractible loop based at $x$ and $L_\gamma$ is
the linear transformation of $\cV_x$ obtained by parallel translation
around $\gamma$.  Just as with infinitesimal holonomy, for the
tractor connection of a conformal manifold we have that if $z\in \cG_x$,
then $\ev_z (\Hol_x^0(\cT,\nabla))$ is an isomorphic copy of  
$\Hol_x^0(\cT,\nabla)$ in $\Aut T_z(\cG\times \R)$. 
\begin{corollary}\label{cor1}
Let $(M,c)$ be an odd-dimensional real-analytic conformal
manifold and $\gt$ a real-analytic ambient metric for $(M,c)$.   
If $x\in M$ and $z\in \cG_x$, then
$$
\ev_z(\Hol_x^0(\cT,\nabla))= \Hol_z^0(T\cGt,\nt).   
$$
\end{corollary}
\noindent
Corollary~\ref{cor1} follows from Theorem~\ref{main} because
$\ev_z(\Hol_x^0(\cT,\nabla))$ and $\Hol_z^0(T\cGt,\nt)$ are connected Lie 
subgroups of $\Aut T_z\cGt$ with the same Lie algebra 
$\ev_z(\fh_x) = \fht_z$.       

The tractor bundle $\cT$ carries a tractor metric $h$ of signature
$(p+1,q+1)$ which is parallel 
with respect to $\nabla$.  So by choosing a frame
for $\cT_x$, one can identify $\Hol^0(\cT,\nabla)$ with a subgroup of  
$\operatorname{SO}_e(p+1,q+1)$ which is well-defined up to conjugacy 
independently of $x$  
and the choice of frame (assuming $M$ is connected).  Corollary~\ref{cor1}
immediately implies:  

\begin{corollary}\label{cor2}
Let $(M,c)$ be an odd-dimensional connected real-analytic conformal
manifold.  Then its restricted conformal holonomy group
$\Hol^0(\cT,\nabla)\subset \operatorname{SO}_e(p+1,q+1)$ is realizable as
the restricted holonomy group of a real-analytic pseudo-Riemannian manifold
of signature $(p+1,q+1)$.  
\end{corollary}
\noindent
Corollary~\ref{cor2} is interesting because of the wealth of known
information concerning pseudo-Riemannian holonomy (in particular, Berger's
list) and the restriction it places on conformal holonomy groups.  

If a pseudo-Riemannian manifold admits a nonzero parallel tensor field,
then its holonomy group is constrained to lie in the isotropy group
consisting of the linear transformations preserving the   
tensor at a point.  Of course, many interesting pseudo-Riemannain holonomy
groups arise in this fashion.  Likewise, interesting classes of
conformal manifolds are characterized by admitting a parallel
tractor-tensor field (i.e. a section of $\otimes^r\cT^*$ for some $r\geq 1$)
of a particular algebraic type.  A precursor to Theorem~\ref{main} is the
result 
of \cite{GW} asserting that a parallel tractor-tensor field on a conformal
manifold admits an extension to the ambient space which is parallel with
respect to the ambient metric (to infinite order for $n$ odd, to order
$n/2-1$ for $n$ even).  This result was one motivation for our
consideration of the 
question of equality of infinitesimal holonomy in general.

In order to prove Theorem~\ref{main}, one must express the
ambient connection and its curvature in tractor terms.  The paper \cite{CG} 
showed how the tractor bundle and connection could be written in ambient
terms.  This gives the inclusion $\ev_z(\fh_x) \subset \fht_z$ 
in Theorem~\ref{main}.  The paper \cite{GP} 
reversed the direction and 
showed how to express the full ambient curvature and its covariant
derivatives in 
terms of tractor calculus.  Our proof of the reverse inclusion in
Theorem~\ref{main}, i.e. of ambient holonomy in tractor holonomy, is based
on these relations.     

In \S\ref{ambtrac} we review the ambient metric construction and the 
realization of the tractor bundle and connection in ambient terms.  In
\S\ref{holonomy} we discuss infinitesimal holonomy and prove
Theorem~\ref{main}, in the process recalling the tractor expressions for 
the ambient curvature and connection.

\section{Ambient Metrics and Tractors}\label{ambtrac} 

We begin by reviewing background material concerning ambient metrics and 
tractors.  The main reference for the material
on ambient metrics is \cite{FG2}.  References for the ambient
formulation of tractors are \cite{CG} and \cite{GP}. 

Let $(M,c)$ be a conformal manifold of dimension $n\geq 3$ and signature
$(p,q)$, $p+q=n$.  Metrics in the conformal class $c$ are sections of the 
metric bundle $\cG:=\{(x,g_x):x\in M, g\in c\}\subset S^2T^*M$.  
Let $\pi:\cG\rightarrow M$ denote the projection and
$\delta_s:\cG\rightarrow \cG$ the  
dilations defined by $\delta_s(x,g_x)=(x,s^2g_x)$, $s>0$.  Let   
$T=\frac{d}{ds}\delta_s|_{s=1}$ be the infinitesimal generator of the 
dilations.  There is a tautological symmetric 2-tensor  
${\bf g}$ on $\cG$ defined for $X$, $Y\in T_{(x,g_x)}\cG$ by 
${\bf g}(X,Y)=g_x(\pi_*X,\pi_*Y)$. 

Regard $\cG$ as a hypersurface in $\cG\times \R$ via   
$\iota(z)=(z,0)$, $z\in\cG$.  The variable in the $\R$ factor is denoted
$\rho$.   
A straight pre-ambient metric 
for $(M,c)$ is a smooth metric $\gt$ of signature $(p+1,q+1)$ on a 
dilation-invariant neighborhood $\cGt$ of $\cG$ satisfying 
\begin{enumerate}
\item[(1)] $\delta_s^* \gt =s^2 \gt\quad$ for $s>0$;
\item[(2)] $\iota^* \gt={\bf g}$;
\item[(3)] $\nt T =Id$, where $Id$ denotes the identity endomorphism and 
$\nt$ the Levi-Civita connection of $\gt$.
\end{enumerate}

If $n$ is odd, an ambient metric for $(M,c)$ is a straight pre-ambient
metric for $(M,c)$ such that $\Ric(\gt)$ vanishes to infinite order on
$\cG$.  (The straightness condition (3) is automatic to infinite 
order, but it is convenient to include it in the definition.)  There exists
an ambient metric for $(M,c)$ and it is unique to 
infinite order up to pullback by a diffeomorphism defined on a
dilation-invariant neighborhood of $\cG\times\R$ which commutes with
dilations and which restricts to the identity on $\cG$.    
If $M$ is a real-analytic manifold and there is a real-analytic metric in
the conformal class,
then there exists a real-analytic 
ambient metric for $(M,c)$ satisfying $\Ric(\gt)=0$ on some
dilation-invariant $\cGt$ as above.

In order to formulate the definition of ambient metrics for $n$ even, 
if $S_{IJ}$ is a symmetric 2-tensor field on an 
open neighborhood of $\cG$ in $\cG \times  \R$ and 
$m \geq 0$, we write $S_{IJ} = O^+_{IJ}( \rho^m)$ if
$S_{IJ} = O(\rho^m)$ and 
for each point $z\in \cG$, the symmetric 2-tensor $(\iota^*(\rho^{-m}S))(z)$  
is of the form 
$\pi^*s$ for some symmetric 2-tensor $s$ at $x=\pi(z)\in M$ satisfying 
$\operatorname{tr}_{g_x}s = 0$.  
If $n$ is even, an ambient metric for $(M,c)$ is a straight pre-ambient
metric such that $\Ric(\gt)=O^+_{IJ}( \rho^{n/2-1})$.
There exists an ambient metric for $(M,c)$ and it is unique up to addition 
of a term which is $O^+_{IJ}( \rho^{n/2})$ and up to pullback by a
diffeomorphism defined on a dilation-invariant 
neighborhood of 
$\cG$ which commutes with dilations and which restricts to the 
identity on $\cG$.  
For $n$ even, a conformally invariant tensor, the ambient obstruction  
tensor, obstructs the existence of smooth solutions to
$\Ric(\gt)=O(\rho^{n/2})$.

Let $(M,c)$ be a conformal manifold with metric bundle  
$\cG\stackrel{\pi}{\rightarrow} M$.  For $x\in M$, write
$\cG_x=\pi^{-1}(\{x\})$ for the fiber of $\cG$ over $x$.  
Recall that the bundle $\cD(w)$ of conformal densities of weight $w\in \C$
has fiber 
$\cD_x(w)=\{f:\cG_x\rightarrow \C:(\delta_s)^*f=s^wf,\;s>0$\}.  Thus  
sections of $\cD(w)$ on $M$ are functions on $\cG$ homogeneous of
degree $w$. 

The standard tractor bundle and its normal connection can be similarly
realized in terms of homogeneous vector fields on $\cG_x$.   
As described in the introduction, the standard tractor bundle can be
realized as the rank $n+2$ vector bundle 
$\cT\rightarrow M$ with fiber over $x$ given by \eqref{tracdef}.
If $\gt$ is an ambient metric for $(M,c)$ and if $U$,
$W\in \cT_x$, then $\gt(U,W)$ is 
homogeneous of degree 0 on $\cG_x$, i.e. $\gt(U,W)\in \R$.  
Therefore $h(U,W)=\gt(U,W)$ 
defines a metric $h$ of signature $(p+1,q+1)$ on $\cT$, the tractor metric.
Since $T$ is homogeneous of degree 0 with  
respect to the $\delta_s$, it defines a section of 
$\cT(1)$, where in general we denote the effect of tensoring a bundle with
$\cD(w)$ by appending $(w)$.  The set of $U$ 
in \eqref{tracdef} which at each point of $\cG_x$ is a  
multiple of $T$ determines a subbundle of $\cT$ which we denote  
$\span\{T\}$.  Its orthogonal complement $\span\{T\}^\perp$ is 
the set of $U$ which at each point of $\cG_x$ is tangent 
to $\cG$.  This gives the filtration
\begin{equation}\label{filtration}
0\subset \text{span}\{T\}\subset \span\{T\}^\perp\subset \cT.
\end{equation}

In order to realize the tractor connection, observe that 
$\pi_*:T\cG\rightarrow TM$ induces a realization of the 
tangent bundle $TM$ as
\begin{equation}\label{tangentbundle}
T_xM=\left\{\eb\in\Gamma(T\cG\,\big{|}_{\cG_x}):
(\delta_s)^*\eb=\eb,\;s>0\right\}\Big{/} \span\{T\},  
\end{equation}
where here $\span\{T\}$ really means the constant multiples of $T$. 
If $\eta\in T_xM$, choose $\eb\in\Gamma(T\cG\,\big{|}_{\cG_x})$ representing 
$\eta$.  We will call such an $\eb$ an invariant lift of $\eta$.  Let $\gt$
be an ambient metric for $(M,c)$ and $\nt$ its Levi-Civita connection.  If
$U$ is a section of $\cT$ near $x$, then 
$\nt_{\eb} U\in \Gamma(T\cGt\,\big{|}_{\cG_x})$ makes sense since $U$ is
defined on $\cG$ and $\eb$ is tangent to $\cG$.  
The straightness of $\gt$ and the homogeneity of $U$ imply that 
$\nt_TU=0$.  Therefore  
$\nt_{\eb} U$ is independent of the choice of invariant lift $\eb$.  Also   
$\nt_{\eb} U$ has the same homogeneity as $U$, so 
$\nt_{\eb} U$ defines an element of $\cT_x$.  This realizes the tractor  
connection $\nabla$ on $\cT$:  
\begin{equation}\label{connection}
\nabla_\eta U = \nt_{\eb} U.    
\end{equation}
The tractor metric $h$ is parallel with respect to $\nabla$ since    
$\nt\gt=0$.  
These realizations of the tractor metric and connection depend on  
the choice of ambient metric $\gt$.  But the realizations obtained by
changing $\gt$ by a diffeomorphism are equivalent.   

The realization \eqref{tracdef} of the tractor bundle induces the following
realizations of the bundles of cotractor-tensors: 
\begin{equation}\label{tensors}
(\otimes^r\cT^*)_x=\left\{\chi\in\Gamma\big(\otimes^rT^*\cGt\,\big{|}_{\cG_x}\big):   
(\delta_s)^*\chi=s^r\chi,\;s>0\right\},\quad r\in \N
\end{equation}
as well as the realization \eqref{endos} of the bundle of tractor
endomorphisms.   
The induced tractor connections on these bundles are also given in terms of
the ambient connection and an invariant lift $\eb$ as in
\eqref{connection}.     
Throughout this paper we will identify weighted tractor-tensors with
homogeneous sections of bundles on $\cG$ as in \eqref{tracdef},
\eqref{endos}, \eqref{tensors}.  

The curvature $R$ of the tractor connection can be expressed in terms of  
the curvature $\Rt$ of an ambient metric.  We have 
$R:\Lambda^2TM\rightarrow \End \cT$ and 
$\Rt:\Lambda^2T\cGt\rightarrow \End T\cGt$.  
It is a fact that $T\into \Rt =0$ on $\cG$, where the contraction can be 
taken in any of the three lower indices.  So if $\eta_1$, $\eta_2\in T_xM$
and $\eb_1$, $\eb_2\in \Gamma\big(T\cG|_{\cG_x}\big)$ are invariant lifts,
then  
$\Rt(\eb_1,\eb_2)\in \Gamma\big(\End T\cGt|_{\cG_x}\big)$ is independent of
the 
choices of $\eb_1$, $\eb_2$.  Moreover, $\Rt(\eb_1,\eb_2)$ is homogeneous
of degree $0$ with respect 
to the $\delta_s$, so it realizes an element of $\End\cT_x$, and one has 
\begin{equation}\label{curvature}
R(\eta_1,\eta_2)=\Rt(\eb_1,\eb_2).   
\end{equation}

We follow usual notational conventions.  We label tensors on the ambient
space and therefore also tractors with capital Latin indices and vectors on  
$M$ with lower case Latin indices.  We use $\cE$ to denote the space of
smooth sections of a bundle on $M$, the bundle specified by
the accompanying indices.  Just as with the bundles themselves, we denote
the spaces of sections of the corresponding weighted bundles by appending 
$(w)$.  The notation $\cE^\Phi(w)$ signifies the space of sections of a
generic weighted tractor bundle, where $\Phi$ denotes an arbitrary
collection of upper and lower capital 
indices.  If $\Phi$ consists of $r$ upper indices and $s$ lower indices, we
denote by $\cEt^\Phi(w)$ the space of sections of $(\otimes^rT\cGt)\otimes
(\otimes^sT^*\cGt)$ on $\cGt$ of the same homogeneity degree as sections of   
$\cE^\Phi(w)$, i.e. of homogeneity degree $w-r+s$.
Ambient/tractor indices are raised and lowered using the ambient/tractor
metric $\gt_{AB}$/$h_{AB}$ and lower case indices using the conformal
metric ${\bf g}_{ij}\in \cE_{ij}(2)$.    

A choice of metric $g$ in the conformal class induces a splitting of the
cotractor bundle
\begin{equation}\label{split}
\cT^*=\cD(-1)\oplus T^*M(1)\oplus \cD(1).
\end{equation}
This is the formulation in the original definition of the tractor bundle in 
\cite{BEG}.  It can also be viewed in terms of the ambient realization 
by putting $\gt$ in normal form relative to $g$ (see \cite{GP} or
\cite{GW}).  The three inclusions determined by this splitting determine
sections  
$$
X_A\in \cE_A(1),\qquad Z_A{}^i\in \cE_A{}^i(-1),\qquad Y_A\in \cE_A(-1)
$$
so that
\begin{equation}\label{splitting}
V_A=\varphi X_A + \psi_i Z_A{}^i + \rho Y_A 
\end{equation}
corresponds to $V_A = (\varphi,\psi_i,\rho)\in \cE(-1)\oplus \cE_i(1)\oplus   
\cE(1)$.  The sections $Y_A$ and $Z_A{}^i$ are scale-dependent, i.e. they
depend on the choice of $g$, while $X_A$ is scale-independent:  
$X^A\in \cE^A(1)$ is   
another notation for the weighted tractor defined by the vector field
$T|_{\cG}$.

\section{Holonomy}\label{holonomy}

Recall from the introduction that the infinitesimal holonomy $\fh_x$ of a
vector bundle with connection $(\cV,\nabla)$ on a manifold $M$ is defined
by \eqref{holdef}.  It is useful to consider the corresponding object  
consisting of global sections.  For $k\geq 2$, we define 
\begin{equation}\label{global}
\fh^k_M=\span_{C^\infty(M)}\left\{\nabla_{\eta_l}\nabla_{\eta_{l-1}}\cdots
\nabla_{\eta_3}\big(R(\eta_1,\eta_2)\big): 2\leq l\leq k, \,\,
\eta_1,\ldots, \eta_l\in \fX(M) \right\}
\end{equation}
and
$$
\fh_M=\bigcup_{k\geq 2}\fh^k_M
$$
so that $\fh^k_M$, $\fh_M\subset \Gamma(\End \cV)$.  Clearly  
$\fh_x=\{E(x):E\in \fh_M\}$.  One has   
\begin{equation}\label{comm}
[\fh_M^k,\fh_M^l]\subset \fh_M^{k+l}.
\end{equation}
In fact, the proof in \cite{KN} that $\fh_x$ is a
subalgebra of $\End \cV_x$ establishes the analog of 
\eqref{comm} in the principal bundle setting.

There is an alternate characterization of these spaces 
in terms 
of iterated covariant derivatives with respect to a coupled connection.  If
we choose 
arbitrarily a connection on $TM$ and denote also by $\nabla$ the coupled
connection on $\cV\otimes TM$, then the Leibniz formula and induction
show that  
\begin{equation}\label{holcoupled}
\fh^k_M=\span_{C^\infty(M)}\left\{(\nabla^{l-2}R)(\eta_1,\eta_2,\ldots,\eta_l)
:2\leq l\leq k,\,\,\eta_1,\ldots,\eta_l\in \fX(M)\right\}.     
\end{equation}
$R$ again denotes the curvature of the connection on $\cV$.  Here it is 
viewed as a section of $\Lambda^2T^*M\otimes \End\cV $ and  
$\nabla^{l-2}R$ denotes its iterated covariant derivative with respect to   
the coupled connection.  

If $(M,c)$ is a conformal manifold, we take $\cV=\cT$ to be the
tractor bundle with its normal connection and we denote the corresponding  
spaces by $\fh^k_M$, $\fh_M$.  As usual, via our realization \eqref{endos}
we identify elements   
of $\fh_M$ as global sections of $\End T(\cG\times \R)\,\big{|}_{\cG}$
which are homogeneous of degree 0 with respect to the $\delta_s$.   

For the ambient metric we modify the definition slightly to respect
homogeneity.  If $n$ is odd and $\gt$ is an ambient metric for $(M,c)$, we
define for $k\geq 2$
\begin{equation}\label{ambientglobal}
\fht{}^k_M=\span_{C^\infty(M)}\left\{\nt_{\xt_l}\nt_{\xt_{l-1}}\cdots
\nt_{\xt_3}\big(\Rt(\xt_1,\xt_2)\big)\big|_\cG: 2\leq l\leq k, \,\,  
\xt_1,\ldots, \xt_l\in \fX_0(\cGt) \right\},
\end{equation}
where $\fX_0(\cGt)$ denotes the space of smooth vector fields on $\cGt$
homogeneous of degree 0 with respect to the $\delta_s$ and $C^\infty(M)$ is
viewed as the subspace of $C^\infty(\cG)$ of functions homogeneous of
degree 0.  Observe that by definition, 
$\fht{}^k_M\subset \Gamma\big(\End T\cGt\,\big{|}_{\cG}\big)$ consists
of sections which are homogeneous of degree 0.  If $n$ is even, we again
define $\fht{}^k_M$ for $k\geq 2$ by \eqref{ambientglobal}, except that we
require that at most $n/2-2$ of the $\xt_i$ are somewhere transverse to
$\cG$.  For general $n$, we then set 
$$
\fht_M=\bigcup_{k\geq 2}\fht{}^k_M.  
$$
As above, $\fht_M$ also has a description in terms of iterated derivatives
of curvature: 
\begin{equation}\label{ambientiterated}
\fht{}^k_M=\span_{C^\infty(M)}\left\{(\nt^{l-2}\Rt)(\xt_1,\xt_2,\ldots,\xt_l)\big|_\cG 
:2\leq l\leq k,\,\,\xt_1,\ldots,\xt_l\in \fX_0(\cGt)\right\},
\end{equation}
where now we take the coupling connection on $T\cGt$ also to be the
Levi-Civita connection $\nt$.  As usual, for $n$ even we require that 
at most $n/2-2$ of the $\xt_i$ are somewhere transverse to $\cG$.  In this
case, the equivalence of the descriptions \eqref{ambientglobal} and  
\eqref{ambientiterated} only holds for $k\leq n/2-1$, since
$\nt_{\xt}\et$ can be transverse to $\cG$ when both $\xt|_{\cG}$ and
$\et|_{\cG}$ are tangent to $\cG$.  

We claim that $\fht_z=\{E(z): E\in \fht_M\}$.  To see this, choose a frame
$\zt_0, \zt_1, \ldots, \zt_{n+1}$ for $T\cGt$ near $z$ such that
$\zt_A|_{\cG}$ is tangent to $\cG$ for $1\leq A\leq n+1$, and such that
each $\zt_A$ is homogeneous of degree $0$ with respect to the $\delta_s$.  
By writing each $\xt_i$ in \eqref{ambhol} as a linear combination of the
$\zt_A$, it is not hard to see that 
$$
\fht_z=\span_\R\left\{\nt_{\zt_{A_k}}\nt_{\zt_{A_{k-1}}}\cdots
\nt_{\zt_{A_3}}\left(\Rt(\zt_{A_1},\zt_{A_2})\right)(z):k\geq 2 \right\},
$$
where for $n$ even at most $n/2-2$ of the indices $A_1,\cdots, A_k$ are
equal to $0$.  It follows immediately that $\fht_z=\{E(z): E\in \fht_M\}$.     

In light of these observations, it is clear that Theorem~\ref{main} is a
consequence of the following theorem.  
\begin{theorem}\label{mainglobal}
Let $(M,c)$ be a conformal manifold of dimension $n\geq 3$ and $\gt$ an 
ambient metric for $(M,c)$.  Then 
$$
\fh_M = \fht_M.  
$$
\end{theorem}

The inclusion $\fh_M \subset \fht_M$ follows immediately from the
ambient realizations of the tractor connection and curvature.  If
$\eta_1,\ldots,\eta_k\in \fX(M)$ and $\eb_1,\ldots,\eb_k$  
are invariant lifts, then \eqref{connection}, \eqref{curvature} give  
\begin{equation}\label{equal}
\nabla_{\eta_k}\nabla_{\eta_{k-1}}\cdots
\nabla_{\eta_3}\left(R(\eta_1,\eta_2)\right)
=\nt_{\eb_k}\nt_{\eb_{k-1}}\cdots
\nt_{\eb_3}\left(\Rt(\eb_1,\eb_2)\right),
\end{equation}
so $\fh_M \subset \fht_M$.  
The right-hand side is in $\fht_M$ also for $n$ even since none of the
$\eb_i$ are transverse to $\cG$.  

We remark that \eqref{equal} is already sufficient to prove
Theorem~\ref{mainglobal}, and therefore also Theorem~\ref{main}, when 
$n=4$.  In fact, when $n=4$, each $\xt_i|_{\cG}$ in \eqref{ambientglobal}
is required to be everywhere 
tangent to $\cG$, so is an invariant lift of some $\eta_i\in \fX(M)$.    

To prove the opposite inclusion $\fht_M \subset \fh_M$, we must rewrite
expressions of the form $\nt_{\xt_l}\nt_{\xt_{l-1}}\cdots
\nt_{\xt_3}\big(\Rt(\xt_1,\xt_2)\big)\big|_\cG$ 
purely in tractor terms when the $\xt_i$ are allowed to be   
transverse to $\cG$.  We do this using tractor representations of the   
curvature and connection of the ambient metric derived in \cite{GP}.  These 
representations are expressed in terms of the splitting \eqref{split},
\eqref{splitting} of the cotractor bundle determined by a choice of metric  
$g\in c$.  Consider first the case $n$ odd.

\bigskip
\noindent
{\it Proof of Theorem~\ref{mainglobal} for $n$ odd}.  
We show by induction on $k\geq 2$ that $\fht{}^k_M\subset \fh_M$.      
For $k=2$, we use the tractor expression for ambient curvature
$$
\Rt_{AB}{}^P{}_Q\big|_\cG 
= Z_A{}^aZ_B{}^bR_{ab}{}^P{}_Q
-\frac{2}{n-4}X_{[A} Z_{B]}{}^b\nabla^cR_{cb}{}^P{}_Q.
$$
This is (13), (35) of \cite{GP}.  The $\nabla^c$ on
the right-hand side refers to the connection obtained by coupling the
tractor connection with the Levi-Civita connection of the chosen
representative metric $g$.  Now $\fht{}^2_M$ is spanned by contractions of
the left-hand side against $\xt_1^A\xt_2^B$, where $\xt_1$, $\xt_2\in
\fX_0(\cGt)$.  It is evident that after such a contraction, the first term
on the right-hand side is in $\fh_M^2$.  For the second term, write 
$\nabla^cR_{cb}{}^P{}_Q={\bf g}^{cd}\nabla_cR_{db}{}^P{}_Q$ and introduce a  
partition of unity subordinate to a covering of $M$ in each open set of
which ${\bf g}^{cd}$ can be expressed as a smooth linear combination of
tensor products of vector fields.  It follows that after contraction with   
$\xt_1^A\xt_2^B$, the second term is in  
$\fh_M^3$.  Thus the initial $k=2$  step of the induction is
established.   

The induction step for higher $k$ will be carried out using the tractor-$D$ 
operator.  If $\Phi$ denotes an arbitrary collection of upper and/or lower
tractor indices, then  
$$
D_A:\cE^\Phi(w)\rightarrow \cE^\Phi{}_A(w-1)
$$ 
is defined in
terms of the splitting determined by a representative metric $g$ by:
\begin{equation}\label{Ddef}
D_AV=w(n+2w-2)Y_AV+(n+2w-2)Z_A{}^a\nabla_aV-X_A\square V,
\end{equation}
where $\square V=\nabla^i\nabla_iV+wJV$ and $J=\frac{R}{2(n-1)}$.  $D_A$
can also be expressed in ambient terms:
\begin{equation}\label{ambientD}
D_AV = (n+2w-2)\nt_A\Vt\big|_{\cG} -X_A(\Dt \Vt)\big|_{\cG}.    
\end{equation}
These are (8), (31) of \cite{GP}.  On the right-hand side,
$\Vt\in\cEt^\Phi(w)$ is an arbitrary homogeneous  
extension of $V\in \cE^\Phi(w)$ and $\Dt$ denotes the ambient Laplacian
acting 
on the corresponding space of tensors:  $\Dt = \nt^I\nt_I$.  The 
expression on the right-hand side turns out to be independent of the choice
of $\Vt$.    

Assume now that $k\geq 2$ and $\fht{}^k_M\subset \fh_M$.  According to  
\eqref{ambientiterated}, in order to prove that 
$\fht{}^{k+1}_M\subset \fh_M$, it suffices to show that 
$(\xt_1^{A}\xt_2^B\cdots \xt_{k+1}^E\nt^{k-1}_{A\cdots
  C}\Rt_{DE}{}^P{}_Q)|_\cG \in \fh_M$ for $\xt_1,\ldots,\xt_{k+1}\in 
\fX_0(\cGt)$.  Set $\xi^A_s=\xt^A_s|_{\cG}\in \cE^A(1)$, $1\leq s\leq
k+1$.  

Define
$$
\Vt=\nt^{k-2}_{B\cdots C}\Rt_{DE}{}^P{}_Q\in
\cEt_{B\cdots E}{}^P{}_Q(-k) 
$$
and rewrite \eqref{ambientD} as  
$$
(n-2k-2)\nt_A\Vt\big|_\cG = D_AV + X_A(\Dt \Vt)\big|_{\cG} 
$$
where $V:=\Vt|_\cG\in \cE_{B\cdots E}{}^P{}_Q(-k)$.  
Since the coefficient $(n-2k-2)$ is nonzero for $n$ odd, it suffices to
show that 
\begin{equation}\label{firstone}
\xi_1^A\cdots \xi_{k+1}^ED_{A}V_{B\cdots E}{}^P{}_Q\in \fh_M   
\end{equation}
and
\begin{equation}\label{secondone}
\xt_2^B\cdots \xt_{k+1}^E\Dt \Vt_{B\cdots E}{}^P{}_Q\big|_{\cG}\in \fh_M.  
\end{equation}

For \eqref{firstone}, contract \eqref{Ddef} against 
$\xi_1^A\cdots \xi_{k+1}^E$.  The first term on the right-hand
side gives a multiple of 
$$
(\xi_1^AY_A)\, \xi_2^B\cdots \xi_{k+1}^EV_{B\cdots E}{}^P{}_Q,  
$$
which is in $\fh_M$ by the induction hypothesis.  The second term on the
right-hand side gives a multiple of 
$$
(\xi_1^AZ_A^a)\xi_2^B\cdots \xi_{k+1}^E\nabla_a
  V_{B\cdots E}{}^P{}_Q.   
$$
If we set $\eta^a=\xi_1^AZ_A{}^a$, then this can be rewritten as 
\[
\begin{split}
\eta^a\xi_2^B\cdots \xi_{k+1}^E&\nabla_a V_{B\cdots E}{}^P{}_Q\\  
&=\nabla_\eta (\xi_2^B\cdots \xi_{k+1}^EV_{B\cdots E}{}^P{}_Q) 
-\sum_{s=2}^{k+1}\xi_2^B\cdots (\nabla_\eta\xi_s^R) \cdots\xi_{k+1}^E 
V_{B\cdots R\cdots E}{}^P{}_Q.
\end{split}
\]
The induction hypothesis shows that 
$\xi_2^B\cdots \xi_{k+1}^EV_{B\cdots E}{}^P{}_Q\in \fh_M$, so we conclude
that 
$\nabla_\eta (\xi_2^B\cdots \xi_{k+1}^EV_{B\cdots E}{}^P{}_Q)\in \fh_M$.   
Each term in the sum on the right-hand side is clearly in $\fh_M$ by the
induction hypothesis.  Thus the contraction of the second term of the
right-hand side of \eqref{Ddef} is in $\fh_M$.  The third term of
\eqref{Ddef} is handled similarly, namely by expanding the difference 
$$
\xi_2^B\cdots \xi_{k+1}^E\nabla^c\nabla_cV_{B\cdots E}{}^P{}_Q
-\nabla^c\nabla_c(\xi_2^B\cdots \xi_{k+1}^EV_{B\cdots E}{}^P{}_Q)
$$
using the Leibniz rule and introducing a partition of unity to 
rewrite sections of tensor product bundles as sums of tensor products of
sections of the factors as in the proof in the case $k=2$.  This
concludes the proof of \eqref{firstone}.    

It remains to prove \eqref{secondone}.  Now $\Dt\Vt=\Dt\nt^{k-2}\Rt$.  It
is well-known that the Laplacian of an iterated covariant derivative of
the curvature tensor of a Ricci-flat metric can be reexpressed  
as a linear combination of quadratic terms in curvature by commuting both
derivatives in $\Dt$ all the way to the right and applying the second
Bianchi identity.  We will argue using the induction hypothesis that each 
commutator term is already in $\fh_M$.  

Write 
$$
\Dt\nt^{k-2}\Rt_{DE}{}^P{}_Q
=\gt^{IJ}\nt_I\nt_J\nt^{k-2}\Rt_{DE}{}^P{}_Q.
$$
First commute $\nt_J$ to the right of all derivatives in $\nt^{k-2}$.
Modulo commutator terms, one obtains
$$
\gt^{IJ}\nt_I\nt^{k-2}\nt_J\Rt_{DE}{}^P{}_Q=
\gt^{IJ}\nt_I\nt^{k-2}\nt_D\Rt_{JE}{}^P{}_Q
+\gt^{IJ}\nt_I\nt^{k-2}\nt_E\Rt_{DJ}{}^P{}_Q.   
$$
Now commuting $\nt_I$ all the way to the right shows that modulo
commutators the above is equal to
$$
\gt^{IJ}\nt^{k-2}\nt_D\nt_I\Rt_{JE}{}^P{}_Q
+\gt^{IJ}\nt^{k-2}\nt_E\nt_I\Rt_{DJ}{}^P{}_Q.
$$
This vanishes on $\cG$ by the second Bianchi identity and the 
infinite-order vanishing of $\Ric{\gt}$.   

To analyze the commutator terms, it is convenient to suppress writing the
$\End T\cGt$ indices ${}^P{}_Q$.  We will denote by $\cRt_{BC}$ the
curvature tensor of $\gt$ viewed as an $\End T\cGt$-valued section of
$\Lambda^2T^*\cGt$.  If $U$ is an $\End T\cGt$-valued section of
$\otimes^rT^*\cGt$ and $V$ is an $\End T\cGt$-valued section of
$\otimes^sT^*\cGt$, we will denote by $[U,V]$ the 
$\End T\cGt$-valued section of $\otimes^{r+s}T^*\cGt$ which is the 
commutator in the $\End T\cGt$ indices and the tensor product in the 
$T^*\cGt$ indices.  The Leibniz formula gives 
\begin{equation}\label{leibnitz}
\nt[U,V]=[\nt U,V]+[U,\nt V].
\end{equation}
The Ricci identity for commuting covariant derivatives can be written  
\begin{equation}\label{ricciid}
[\nt_B,\nt_C]U = \cRt_{BC}.U+[\cRt_{BC},U], 
\end{equation}
where $\cRt_{BC}.U$ denotes the action of the endomorphism
$\cRt_{BC}$ on the $\otimes^rT^*\cGt$ indices of $U$.  

Every commutator which arose in the above argument was of the form
$$
\nt^i[\nt_B,\nt_C]\nt^j\cRt
$$
for some choice of indices $B$, $C$, where $i\geq 0$, $j\geq 0$, 
and $i+j=k-2$.  Express the commutator $[\nt_B,\nt_C]\nt^j\cRt$ using
\eqref{ricciid} with $U=\nt^j\cRt$.  The first term on the right-hand side
of \eqref{ricciid} gives terms of the form $\nt^i(\cRt.\nt^j\cRt)$.
Expanding the $\nt^i$ with the Leibniz rule, it is clear that
one obtains a sum of terms, each of which has the form
\begin{equation}\label{firsttype}
\contr(\nt^p\cRt\otimes\nt^q\cRt)
\end{equation}
with $p\geq 0$, $q\geq 0$, and $p+q=k-2$.  Here $\contr$ indicates a 
single contraction of the upper $\End T^*\cGt$ index of $\nt^p\Rt$ against  
one of the $\otimes^{q+2}T^*\cGt$ indices of $\nt^q\cRt$.  In particular,
the suppressed $\End T^*\cGt$ indices are those on $\nt^q\cRt$.  The
second term on  
the right-hand side of \eqref{ricciid} gives terms of the form 
$\nt^i[\cRt,\nt^j\cRt]$.  Expanding the $\nt^i$ using \eqref{leibnitz}, one
obtains a sum of terms of the form 
\begin{equation}\label{secondtype}
[\nt^p\cRt,\nt^q\cRt],
\end{equation}
again with $p\geq 0$, $q\geq 0$, and $p+q=k-2$.   

We need to show \eqref{secondone}.  Suppressing the $\End T\cGt$ indices,
we have 
$$
\xt_2^B\cdots \xt_{k+1}^E\Dt\Vt_{B\cdots E}|_{\cG}
=\xt_2^B\cdots \xt_{k+1}^E\gt^{IJ}\nt_I\nt_J
\nt^{k-2}_{B\cdots C}\cRt_{DE}|_{\cG}.
$$
Upon commuting $\nt_I$ and $\nt_J$ to the right as described above, it
follows that this may be written as a sum of contractions of terms of the
form 
\eqref{firsttype}, \eqref{secondtype} against $\xt_i$ and $\gt^{IJ}$ with
all indices 
contracted except for the suppressed $\End T\cGt$ indices.  In a term  
\eqref{firsttype}, the free $\End T\cGt$ indices are those on the
second factor $\nt^q\cRt$.  Consequently, we can introduce a partition of
unity and express locally the tensor arising from $\gt^{IJ}$, $\nt^p\cRt$, and
the $\xt_i$ which contracts against the other $q+2$ 
indices of $\nt^q\cRt$ as a sum of tensor products of vector fields.  Since
$q\leq k-2$, it follows by the induction hypothesis that all these terms 
are in $\fh_M$ when restricted to $\cG$.  In a term
\eqref{secondtype}, all the indices except the endomorphism indices
are contracted against $\gt^{IJ}$ and the $\xt_i$.  Again use a partition
of unity and express locally $\gt^{IJ}$ as a sum of tensor products of
vector fields.  Then the induction hypothesis implies that the restriction
to $\cG$ of the contractions against $\nt^p\cRt$ and $\nt^q\cRt$ are 
separately in $\fh_M$.  It follows from \eqref{comm}  that the
commutator is also in $\fh_M$.    
\stopthm

\bigskip
\noindent
{\it Proof of Theorem~\ref{mainglobal} for $n$ even}.  
We have already observed that \eqref{equal} is sufficient to prove the case
$n=4$.  So we assume that $n\geq 6$. 
We next observe that the same argument used for $n$ odd applies  
also when $n$ is even to show $\fht{}^{n/2-1}_M\subset \fh_M$.  In fact,
up to this order the relevant constant $n+2w-2$ in \eqref{ambientD} is
nonzero and the argument only uses $\Ric(\gt)=O(\rho^{n/2-1})$.  

For $n\geq 6$ even, we prove $\fht{}^k_M\subset \fh_M$ by induction on $k$, 
beginning with the 
case $k=n/2-1$.  So assume for some $k\geq n/2-1$ that $\fht{}^k_M\subset
\fh_M$ and we will show that $\fht{}^{k+1}_M\subset \fh_M$.   
According to \eqref{ambientglobal}, we have to show that 
$\nt_{\xt_{k+1}}\nt_{\xt_k}\cdots
\nt_{\xt_3}\big(\Rt(\xt_1,\xt_2)\big)\big|_\cG\in \fh_M$ whenever
$\xt_1,\ldots, \xt_{k+1}\in \fX_0(\cGt)$ and at most $n/2-2$ of the $\xt_i$
are somewhere transverse to $\cG$.  Since $k+1\geq n/2$, at least two of
the $\xt_i$ are everywhere tangent to $\cG$.  If $\xt_{k+1}$ is everywhere
tangent to $\cG$, then its restriction to $\cG$ is the invariant lift of
some $\eta\in \fX(M)$.  In this case \eqref{connection} gives 
$$
\nt_{\xt_{k+1}}\nt_{\xt_k}\cdots
\nt_{\xt_3}\big(\Rt(\xt_1,\xt_2)\big)\big|_\cG
=\nabla_\eta\Big(\nt_{\xt_k}\cdots
\nt_{\xt_3}\big(\Rt(\xt_1,\xt_2)\big)\big|_\cG\Big).  
$$
The induction hypothesis shows that 
$\nt_{\xt_k}\cdots \nt_{\xt_3}\big(\Rt(\xt_1,\xt_2)\big)\big|_\cG\in
\fh_M$, from which it follows that 
$\nabla_\eta\Big(\nt_{\xt_k}\cdots
\nt_{\xt_3}\big(\Rt(\xt_1,\xt_2)\big)\big|_\cG\Big) \in \fh_M$ as desired.

If $\xt_i$ is everywhere tangent to $\cG$ for some $i$, $3\leq i\leq k$,
then we can commute $\nt_{\xt_i}$ all the way to the left and reduce to the  
previous case.  Modulo relabeling the indices, each commutator is
of the form 
\[
\begin{split}
\nt_{\xt_{k+1}}\cdots
\nt_{\xt_{j+1}}&[\nt_{\xt_{j}},\nt_{\xt_{j-1}}]\nt_{\xt_{j-2}}\cdots
\nt_{\xt_3}\Rt(\xt_1,\xt_2)\\
=\nt_{\xt_{k+1}}&\cdots 
\nt_{\xt_{j+1}}\nt_{[\xt_{j},\xt_{j-1}]}\nt_{\xt_{j-2}}\cdots
\nt_{\xt_3}\Rt(\xt_1,\xt_2)\\
&+\nt_{\xt_{k+1}}\cdots
\nt_{\xt_{j+1}}[\Rt(\xt_{j},\xt_{j-1}),\nt_{\xt_{j-2}}\cdots
\nt_{\xt_3}\Rt(\xt_1,\xt_2)].
\end{split}
\]
In the first term on the right-hand side, the number of
differentiations has decreased by 1 without increasing the number of
vector fields somewhere transverse to $\cG$, since the commutator of two
vector fields tangent to $\cG$ is also tangent to $\cG$.  So the
restriction to $\cG$ of this term is in $\fh_M$ by the induction
hypothesis.  In the second term on the right-hand side, expand the
derivatives outside the commutator using the Leibniz rule.  One obtains a
linear combination of commutators of covariant derivatives of curvature
endomorphisms.  The restriction to $\cG$ of each such covariant derivative
itself is in $\fh_M$ by the induction hypothesis.  Equation \eqref{comm}
then shows that the commutator is in $\fh_M$.    

Finally we must consider the possibility that none of
$\xt_3,\ldots,\xt_{k+1}$ is everywhere tangent to $\cG$.  (This can only
happen in the beginning case $k=n/2-1$, but we will not use this.)  It  
must be that $\xt_1$ and $\xt_2$ are everywhere tangent to $\cG$.  In
this case, we apply the second Bianchi identity to write  
\[
\begin{split}
\nt_{\xt_3}\Rt(\xt_1,\xt_2)=&\nt_{\xt_1}\Rt(\xt_3,\xt_2)+\nt_{\xt_2}\Rt(\xt_1,\xt_3)\\ 
+& \Rt(\nt_{\xt_3}\xt_1,\xt_2)
+ \Rt(\xt_1,\nt_{\xt_3}\xt_2)
-\Rt(\nt_{\xt_1}\xt_3,\xt_2)
-\Rt(\xt_1,\nt_{\xt_2}\xt_3) 
\\
+& \Rt(\nt_{\xt_1}\xt_2-\nt_{\xt_2}\xt_1,\xt_3).
\end{split}
\]
The terms on the first line of the right-hand side reduce to the previous
case.  The terms on the second line of the right-hand side reduce to the
induction hypothesis since $\xt_1$ and $\xt_2$ are tangential
and at least one occurs as an argument in each term, so the number of
transversal vector fields does not increase.  The last term also 
reduces to the induction hypothesis since
$\nt_{\xt_1}\xt_2-\nt_{\xt_2}\xt_1=[\xt_1,\xt_2]$ is tangential.  
\stopthm


\begin{thebibliography}{BEGM}

\bibitem[A]{Armstrong} S. Armstrong, {\it Definite signature conformal
  holonomy: a complete classification}, {J.\ Geom.\ Phys.} {\bf
  57} (2007), 2024--2048, {\tt arxiv:math/0503388}.  


\bibitem[BEG]{BEG} T. N. Bailey, M. G. Eastwood and A. R. Gover, {\it 
Thomas's structure bundle for conformal, projective and related
structures}, Rocky Mountain J. Math. {\bf 24} (1994), 1191--1217. 


\bibitem[\v{C}G]{CG} A. \v{C}ap and A. R. Gover, {\it Standard
  tractors and the conformal ambient metric construction}, Ann. Global
  Anal. Geom. {\bf 24} (2003), 231--259, {\tt arXiv:math/0207016}.   


\bibitem[FG1]{FG1} C. Fefferman and C. R. Graham, {\it
Conformal invariants,} in {\it The Mathematical Heritage of \'Elie Cartan 
(Lyon, 1984)},  Ast\'erisque, 1985, Numero Hors Serie, 95--116.  

\bibitem[FG2]{FG2} C. Fefferman and C. R. Graham, {\it
The Ambient Metric}, Princeton University Press, 2012, 
{\tt arXiv:math/0710.0919}.

\bibitem[GP]{GP} A. R. Gover and L. J. Peterson, {\it Conformally invariant
  powers of the Laplacian, $Q$-curvature, and tractor calculus},
  Comm. Math. Phys. {\bf 235} (2003), 339--378, {\tt arXiv:math-ph/0201030}.

\bibitem[GW]{GW} C. R. Graham and T. Willse, {\it Parallel tractor
  extension and ambient metrics of holonomy split $G_2$}, J. Diff. Geom.,
  {\bf 92} (2012), 463--505, {\tt arXiv:1109.3504}.  


\bibitem[KN]{KN} S. Kobayashi and K. Nomizu, {\it Foundations of
  Differential Geometry, I}, Interscience Publishers, 1969.  

\bibitem[L]{Leitner} F. Leitner, {\it Conformal Killing forms with
  normalisation condition}, Rend. Circ. Mat. Palermo (2) Suppl. No. 75
  (2005), 279--292. 



\end{thebibliography}
\end{document}